\documentclass[12pt]{article}
\usepackage{amssymb}
\usepackage{times}
\usepackage[T1]{fontenc}
\usepackage{mathrsfs}
\usepackage{latexsym}
\usepackage[dvips]{graphics}
\usepackage{epsfig}
\usepackage{amsmath,amsfonts,amsthm,amscd}
\usepackage{booktabs}
\usepackage[all,cmtip]{xy}

\begin{document}

\theoremstyle{plain} \newtheorem*{smoothc}{Theorem 1}
\theoremstyle{plain} \newtheorem*{corsmoothc}{Corollary}
\theoremstyle{plain} \newtheorem*{powerfullc}{Theorem 2}
\theoremstyle{plain} \newtheorem*{GS2.1}{Theorem A}
\theoremstyle{plain} \newtheorem*{GS2.2}{Theorem B}
\theoremstyle{plain} \newtheorem*{Dchi1}{Lemma 3}
\theoremstyle{plain} \newtheorem*{smoothLfn}{Lemma 4}
\theoremstyle{plain} \newtheorem*{powerfullLfn}{Lemma 5}
\theoremstyle{plain}
\newtheorem{thm}{Theorem}[section] \theoremstyle{plain}
\newtheorem{lem}[thm]{Lemma} \theoremstyle{plain}
\newtheorem{cor}[thm]{Corollary} \theoremstyle{definition}
\newtheorem{defn}[thm]{Definition} \theoremstyle{plain}
\newtheorem{prop}[thm]{Proposition} \theoremstyle{plain}
\newtheorem{exa}[thm]{Example} \theoremstyle{plain}
\newtheorem{exe}{Problem} \theoremstyle{plain}
\newtheorem{conj}[thm]{Conjecture} \theoremstyle{plain}
\newtheorem{cla}[thm]{Claim} \theoremstyle{plain}
\newtheorem{rmk}[thm]{Remark}

\newcommand{\dx}{\text{ } d}
\newcommand{\ZqZ}{\mathbb{Z}/q\mathbb{Z}}
\newcommand{\ZqZstar}{\left(\mathbb{Z}/q\mathbb{Z}\right)^*}
\newcommand{\chimodq}{$\chi$ (mod $q$) }
\newcommand{\bthm}{\begin{thm}}
\newcommand{\ethm}{\end{thm}}
\newcommand{\bcor}{\begin{cor}}
\newcommand{\ecor}{\end{cor}}
\newcommand{\blem}{\begin{lem}}
\newcommand{\elem}{\end{lem}}
\newcommand{\bprop}{\begin{prop}}
\newcommand{\eprop}{\end{prop}}
\newcommand{\bexa}{\begin{exa}}
\newcommand{\eexa}{\end{exa}}
\newcommand{\bdefn}{\begin{defn}}
\newcommand{\edefn}{\end{defn}}
\newcommand{\bexe}{\begin{exe}}
\newcommand{\eexe}{\end{exe}}
\newcommand{\bconj}{\begin{conj}}
\newcommand{\econj}{\end{conj}}
\newcommand{\bcla}{\begin{cla}}
\newcommand{\ecla}{\end{cla}}
\newcommand{\brmk}{\begin{rmk} $\\$}
\newcommand{\ermk}{\end{rmk}}
\newcommand{\bleft}{\begin{flushleft}}
\newcommand{\eleft}{\end{flushleft}}
\newcommand{\bproof}{\begin{proof} $\\$}
\newcommand{\eproof}{$\\$ \end{proof}}

\title{Character sums to smooth moduli are small}
\author{Leo Goldmakher}
\date{}

\maketitle


\begin{abstract}
\noindent Recently \cite{GS1}, Granville and Soundararajan have made
fundamental breakthroughs in the study of character sums. Building
on their work and using estimates on short character sums developed
by Graham-Ringrose \cite{G-R} and Iwaniec \cite{Iw}, we improve the
P\'{o}lya-Vinogradov inequality for characters with smooth conductor.
\end{abstract}


\section{Introduction}

Introduced by Dirichlet to prove his celebrated theorem on primes in
arithmetic progressions (see \cite{DMNT}), Dirichlet characters have
proved to be a fundamental tool in number theory. In particular,
character sums of the form
$$
S_\chi(x) := \sum_{n \leq x} \chi(n)
$$
(where \chimodq is a Dirichlet character) arise naturally in many
classical problems of analytic number theory, from estimating the
least quadratic nonresidue (mod $p$) to bounding $L$-functions.
Recall that for any character \chimodq\!, $\left| S_\chi(x) \right|$
is trivially bounded above by $\varphi(q)$. A folklore conjecture
(which is a consequence of the Generalized Riemann Hypothesis)
predicts that for non-principal characters the true bound should look
like\footnote{Here and throughout we use Vinogradov's notation $f
\ll g$ to mean $f = O(g)$, with variables in subscript to
\indent\indent $\!$indicate dependence of the implicit constant.\newline\newline 
\indent\indent $\!$AMS subject classification: 11L40, 11M06.}
$$
\left| S_\chi(x) \right| \ll_\epsilon \sqrt{x} \cdot q^\epsilon .
$$
Although we are currently very far from
being able to prove such a statement, there have been some
significant improvements over the trivial estimate. The first such
is due (independently) to P\'{o}lya and Vinogradov: they proved that
$$
\left|S_\chi(x)\right| \ll \sqrt{q} \, \log q
$$
(see \cite{DMNT}, pages 135-137). Almost 60 years later, Montgomery
and Vaughan \cite{MV1} showed that conditionally on the Generalized
Riemann Hypothesis (GRH) one can improve P\'{o}lya-Vinogradov to
$$
\left| S_\chi(x) \right| \ll \sqrt{q} \, \log \log q .
$$
This is a best possible result, since in 1932 Paley \cite{Paley} had
given an unconditional construction of an infinite class of
quadratic characters for which the magnitude of the character sum
could be made $\gg \sqrt{q} \, \log \log q$. \newline

In their recent work \cite{GS1}, Granville and Soundararajan give a
characterization of when a character sum can be large; from this
they are able to deduce a number of new results, including
an improvement of P\'{o}lya-Vinogradov (unconditionally) and of
Montgomery-Vaughan (on GRH) for characters of small odd order. In
the present paper we explore a different application of their
characterization. Recall that a positive integer $N$ is said to be
{\em smooth} if its prime factors are all small relative to $N$; if
in addition the product of all its prime factors is small, $N$ is
{\em powerful}. Building on the work of Granville-Soundararajan and
using a striking estimate developed by Graham and Ringrose, we will
obtain (in Section \ref{upperbound}) the following improvement of
P\'{o}lya-Vinogradov for characters of smooth conductor:
\begin{smoothc}
Given \chimodq a primitive character, with $q$ squarefree. For any integer
$n$, denote its largest prime factor by $\mathcal{P}(n)$. Then
$$
\left| S_\chi(x) \right| \ll
\sqrt{q} \, (\log q) \left(\left(\frac{\log \log \log q}{\log \log q}\right)^\frac{1}{2} +
\left(\frac{(\log \log \log q)^2 \log \big(\mathcal{P}(q) \, d(q)\big)}{\log q}\right)^{1/4} \right)
$$
where $d(q)$ is the number of divisors of $q$, and the implied constant is absolute.
\end{smoothc}
From the well-known upper bound $\log d(q) < \frac{\log q}{\log\log
q}$ (see, for example, Ex. 1.3.3 of \cite{MRM}), we immediately
deduce the following weaker but more palatable bound:
\begin{corsmoothc}
Given \chimodq primitive, with $q$ squarefree. Then
$$
\left| S_\chi(x) \right| \ll \sqrt{q} \, (\log q)
\left(
\frac{(\log \log \log q)^2}{\log \log q} + 
\frac{(\log \log \log q)^2 \log \mathcal{P}(q)}{\log q} \right)^\frac{1}{4}
$$
where the implied constant is absolute.
\end{corsmoothc}

For characters with powerful conductor, we can do better by
appealing to work of Iwaniec \cite{Iw}. We prove:
\begin{powerfullc}
Given \chimodq a primitive Dirichlet character with $q$ large and
$$
\text{rad} (q) \leq \exp\left(\left(\log q\right)^{3/4}\right) ,
$$
where the radical of $q$ is defined
$$
\text{rad} (q) := \prod_{p \mid q} p .
$$
Then
$$
\left| S_\chi (x) \right| \ll_\epsilon \sqrt{q} \, (\log q)^{7/8 +
\epsilon} .
$$
\end{powerfullc}

The key ingredient in the proofs of Theorems 1 and 2 is also at the
heart of \cite{GS1}. In that paper, Granville and Soundararajan
introduce a notion of `distance' on the set of characters, and then
show that $\left|S_\chi(x)\right|$ is large if and only if $\chi$ is
close (with respect to their distance) to a primitive character of
small conductor and opposite parity (ideas along these lines had
been earlier approached by Hildebrand in \cite{Hild}, and - in the
context of mean values of arithmetic functions - by Hal\'{a}sz in
\cite{H1, H2}). More precisely, given characters $\chi,\psi$, let
$$
\mathbb{D}(\chi,\psi; y) := \left(\sum_{p \leq y} \frac{1 - \text{Re
} \chi(p) \overline{\psi} (p)}{p} \right)^{\frac{1}{2}} .
$$
Although one can easily furnish characters $\chi \neq \psi$ and a $y$ for which
$\mathbb{D}(\chi,\psi;y) = 0$, all the other properties of a distance function are satisfied;
in particular, a triangle inequality holds:
$$
\mathbb{D}(\chi_1,\psi_1;y) + \mathbb{D}(\chi_2,\psi_2;y) \geq
\mathbb{D}(\chi_1\chi_2,\psi_1\psi_2;y)
$$
(see \cite{GS2} for a more general form of this `distance' and its
role in number theory). Granville and Soundararajan's
characterization of large character sums comes in the form of the
following two theorems:
\begin{GS2.1}[\cite{GS1}, Theorem 2.1]
Given \chimodq primitive, let $\xi \, (\text{mod } m)$ be any
primitive character of conductor less than $(\log q)^{\frac{1}{3}}$
which minimizes the quantity $\mathbb{D}(\chi,\xi;q)$. Then
$$
\left|S_\chi(x)\right| \ll \left(1 - \chi(-1) \xi(-1) \right)
\frac{\sqrt{m}}{\varphi(m)} \, \sqrt{q} \, \log q \,
\exp\left(-\frac{1}{2} \, \mathbb{D}(\chi,\xi;q)^2\right) + \sqrt{q}
\, \left(\log q\right)^{\frac{6}{7}} .
$$
\label{thm2.1}
\end{GS2.1}
\begin{GS2.2}[\cite{GS1}, Theorem 2.2]
Given \chimodq a primitive character, let $\xi \, (\text{mod } m)$
be any primitive character of opposite parity. Then
$$
\max_x \left|S_\chi(x)\right| + \frac{\sqrt{m}}{\varphi(m)} \,
\sqrt{q} \, \log\log q \gg \frac{\sqrt{m}}{\varphi(m)} \, \sqrt{q}
\, \log q \, \exp\left(-\mathbb{D}(\chi,\xi;q)^2\right) .
$$
\label{thm2.2}
\end{GS2.2}
Roughly, the first theorem says that $\left|S_\chi(x)\right|$ is
small (i.e. $\ll \sqrt{q} \, (\log q)^{6/7}$) unless there exists a
primitive character $\xi$ of small conductor and opposite parity,
whose distance from $\chi$ is small (i.e. $\mathbb{D}(\chi,\xi;q)^2
\leq \frac{2}{7} \, \log \log q$); the second theorem says that if
there exists a primitive character $\xi \, (\text{mod } m)$ of small
conductor
and of opposite parity, whose distance from $\chi$ is small,
then $\left|S_\chi(x)\right|$ gets large.
In particular, to improve P\'{o}lya-Vinogradov for a primitive
character \chimodq it suffices (by Theorem A) to find a lower bound
on the distance from $\chi$ to primitive characters of small
conductor and opposite parity. For example, if one can find a
positive constant $\delta$, independent of $q$, for which
\begin{equation}
\mathbb{D}(\chi,\xi;q)^2 \geq \left(\delta + o(1)\right) \log \log q
\label{deltag}
\end{equation}
then Theorem A would immediately yield an improvement of
P\'{o}lya-Vinogradov:
$$
\max_x \left|S_\chi(x)\right| \ll \sqrt{q} \, \left(\log q\right)^{1
- \frac{\delta}{2} + o(1)} .
$$
As it turns out (see Lemma 3.2 of \cite{GS1}), it is not too
difficult to show (\ref{deltag}) holds for $\chi$ a character of odd
order $g$, with $\delta = \delta_g = 1 - \frac{g}{\pi} \, \sin
\frac{\pi}{g}$.
\newline

Thus, to derive bounds on character sums from Theorem A, one must
understand the magnitude of $\mathbb{D}(\chi,\xi;q)$; this is the
problem we take up in Section \ref{lowerbound}. Since
$\mathbb{D}(\chi,\xi;q) = \mathbb{D}(\chi \, \overline{\xi},1;q)$,
we are naturally led to study lower bounds on distances of the form
$\mathbb{D}(\chi, 1; y)$, for $\chi$ a primitive character and $y$ a
parameter with some flexibility. By definition,
$$
\mathbb{D}(\chi,1;y)^2  = \sum_{p \leq y} \frac{1}{p} - \text{Re }
\sum_{p \leq y} \frac{\chi(p)}{p} .
$$
The first sum on the right hand side is well-approximated by
$\log\log y$ (a classical estimate due to Mertens, see pages 56-57
of \cite{DMNT}); we will show that the second sum is comparable to
$\left|L(s_y,\, \chi)\right|$, where
$$
s_y := 1 + \frac{1}{\log y} \, .
$$
To be precise, in Section \ref{lowerbound} we prove:
\begin{Dchi1}
\label{thmDchi1}
For all $y \geq 2$,
$$
\mathbb{D}(\chi,1;y)^2 = \log \left| \frac{\log y}{L(s_y,\chi)}
\right| + O(1) .
$$
\end{Dchi1}
Our problem is now reduced to finding upper bounds on $\left|
L(s,\chi) \right|$ for $s$ slightly larger than 1. This is a
classical subject, and many bounds are available. Thanks to the
remarkable work of Graham and Ringrose \cite{G-R} on short
character sums, a particularly strong upper bound on $L$-functions
is known when the character has smooth modulus; from a slight
generalization of their result we will deduce (in Section
\ref{smoothL}) the following:
\begin{smoothLfn}
Given a primitive character $\chi$ (mod $Q$), let $r$ be any
positive number such that for all $p \geq r$, $\text{ord}_p \, Q
\leq 1$. Let
$$
q' = q_r' := \prod_{p < r} p^{\text{ord}_p \, Q}
$$
and denote by $\mathcal{P}(Q)$ the largest prime factor of $Q$. Then for all $y > 3$,
$$
\left| L(s_y, \chi) \right| \ll \log q' + \frac{\log Q}{\log \log Q} +
\sqrt{(\log Q)\big(\log \mathcal{P}(Q) + \log d(Q)\big)}
$$
where the implied constant is absolute.
\end{smoothLfn}
Using the bound $\log d(q) < \frac{\log q}{\log \log q}$
one deduces the friendlier but weaker bound
$$
\left| L(s_y, \chi) \right| \ll \log q' + \frac{\log Q}{(\log \log
Q)^{1/2}} + \sqrt{(\log Q)\big(\log \mathcal{P}(Q)\big)} .
$$

Lemma 4 will enable us to prove Theorem 1. For the proof of Theorem
2, we need a corresponding bound for $L(s_y,\, \chi)$ when the
conductor of $\chi$ is powerful. In Section \ref{powerfullL}, using
a potent estimate of Iwaniec \cite{Iw} we will prove:
\begin{powerfullLfn}
Given $\chi$ (mod $Q$) a primitive Dirichlet character with $Q$
large and
$$
\text{rad}(Q) \leq \exp \left(2 \left(\log Q\right)^{3/4}\right) .
$$
Then for all $y > 3$,
$$
\left| L(s_y,\, \chi) \right| \ll_\epsilon \left(\log Q\right)^{3/4
+ \epsilon} .
$$
\end{powerfullLfn}

In the final section of the paper, we synthesize our results and prove
Theorems 1 and 2.\newline

{\em Acknowledgements:} I am indebted to Professor Soundararajan for
suggesting the problem in the first place, for encouraging me
throughout, and for making innumerable improvements to my
exposition. I am also grateful to the referee for meticulously
reading the manuscript and catching an important error in the proofs 
of Theorems 1 and 2, to Denis Trotabas and Bob Hough for some helpful 
discussions, and to the Stanford Mathematics Department, where the bulk 
of this project was completed.

\section{The size of $\mathbb{D}(\chi,1;y)$} \label{lowerbound}

How large should one expect $\mathbb{D}(\chi,1;y)$ to be? Before proving Lemma 3
we gain intuition by exploring what can be deduced from GRH.

\begin{prop}
Assume GRH. For any non-principal character $\chi$ (mod $Q$) we have
$$
\mathbb{D}(\chi,1;y)^2 = \log\log y + O(\log\log\log Q) .
$$
\label{condprop}
\end{prop}
\begin{proof}
Since
$$
\sum_{p \leq y} \frac{1}{p} = \log \log y + O(1)
$$
by Mertens' well-known estimate, we need only show that
$$\sum_{p \leq y} \frac{\chi(p)}{p} = O\left(\log \log \log Q\right).$$
We may assume that $y > (\log Q)^6$, else the estimate is trivial.
Recall that on GRH, for all $x > (\log Q)^6$ we have:
$$
\theta(x,\chi) := \sum_{p \leq x} \chi(p) \log p
\ll \sqrt{x} \left(\log Qx\right)^2 \ll x^{5/6}
$$
(such a bound may be deduced from the first formula appearing on
page 125 of \cite{DMNT}). Partial summation now gives
$$
\sum_{(\log Q)^6 < p \leq y} \frac{\chi(p)}{p} =
\int_{(\log Q)^6}^y \frac{1}{t \log t} \dx \theta(t,\chi) \ll
\frac{1}{\log Q} \ll 1
$$
and the proposition follows.
\end{proof}

We now return to unconditional results. Recall that the prime number
theorem gives $\theta(x) := \sum_{p \leq x} \log p \sim x$.
\begin{proof}[Proof of Lemma 3:]
As before, by Mertens' estimate it suffices to show that
\begin{equation}
\text{Re } \sum_{p \leq y} \frac{\chi(p)}{p} =
\log \left|L(s_y, \chi)\right| + O(1)
\label{lem3step1}
\end{equation}
where $s_y := 1 + (\log y)^{-1}$. From the Euler product we know
$$
\log \left|L(s_y, \chi)\right| =
\text{Re } \sum_p \sum_{k=1}^\infty \frac{\chi(p)^k}{k \, p^{k s_y}} =
\text{Re } \sum_p \frac{\chi(p)}{p^{s_y}} + O(1)
$$
so that (\ref{lem3step1}) would follow from
$$
\sum_{p \leq y} \left(\frac{1}{p} - \frac{1}{p^{s_y}} \right) +
\sum_{p > y} \frac{1}{p^{s_y}} \ll 1 .
$$
The first term above is
$$
\sum_{p \leq y} \left(\frac{1}{p} - \frac{1}{p^{s_y}} \right) =
\sum_{p \leq y} \frac{1 - \exp\left(-\frac{\log p}{\log y}\right)}{p} \leq
\frac{1}{\log y} \, \sum_{p \leq y} \frac{\log p}{p} =
\frac{1}{\log y} \, \int_1^y \frac{1}{t} \dx \theta(t) \ll 1
$$
by partial summation and the prime number theorem. A second application of partial summation and
the prime number theorem yields
$$
\sum_{p > y} \frac{1}{p^{s_y}} = \int_y^\infty \frac{1}{t^{s_y} \log
t} \dx \theta(t) \ll 1 .
$$
The lemma follows.
\end{proof}

For a clearer picture of where we are heading, we work out a simple
consequence of this result. Let \chimodq and $\xi$ (mod $m$) be as
in Theorem A. By Lemma 3,
$$
\mathbb{D}(\chi,\xi;q)^2 = \mathbb{D}(\chi\overline{\xi},1;q)^2 =
\log \left| \frac{\log q}{L(s_q,\chi \, \overline{\xi})} \right| +
O(1) \, ,
$$
and Theorem A immediately yields:
\bprop
Let \chimodq be a primitive character, and $\xi$ a character as in Theorem A.
Then 
$$
\left| S_\chi(x) \right| 
\ll
\sqrt{q} \, \sqrt{(\log q) \, \left|L(s_q,\chi \,
\overline{\xi})\right|} + \sqrt{q} \, (\log q)^{6/7} .
$$
\label{TheBound}
\eprop

Thus to improve P\'{o}lya-Vinogradov it suffices to prove
$$
L(s_q, \chi \, \overline{\xi}) = o(\log q) .
$$
This is the problem we explore in the next two sections.

\section{Proof of Lemma 4} \label{smoothL}

We ultimately wish to bound $|L(s_q, \chi \, \overline{\xi})|$; in
this section we explore the more general quantity $|L(s_y, \chi)|$,
where throughout $y$ will be assumed to be at least 3, and $Q$ will
denote the conductor of $\chi$. \newline

By partial summation (see (8) on page 33 of \cite{DMNT}),
$$
L(s_y,\chi) = s_y \int_1^\infty \frac{1}{t^{s_y + 1}} \left(\sum_{n
\leq t} \chi(n)\right) \dx t .
$$
When $t > Q$, the character sum is trivially bounded by $Q$, so that this portion of the integral
contributes an amount $\ll 1$. For $t \leq T$ (a suitable parameter to be chosen later),
we may bound our character sum by $t$, and therefore this portion of the integral contributes an
amount $\ll \log T$. Thus,
\begin{equation}
|L(s_y, \chi)| \ll \left| \int_T^Q \frac{1}{t^2} \left(\sum_{n \leq
t} \chi(n)\right) \dx t \right| + 1 + \log T . \label{LsChi}
\end{equation}
To bound the character sum in this range, we invoke a powerful estimate of Graham and Ringrose
\cite{G-R}. For technical reasons, we need a slight generalization of their theorem:

\begin{thm}[compare to Lemma 5.4 of \cite{G-R}]
Given a primitive character $\chi$ (mod $Q$), with $q'$ and $\mathcal{P}(Q)$
defined as in Lemma 4. Then for any $k \in \mathbb{N}$, writing $K := 2^k$, we have
$$
\left| \sum_{M < n \leq M + N} \chi(n) \right| \ll
N^{1 - \frac{k+3}{8K-2}} \,
\mathcal{P}(Q)^{\frac{k^2 + 3k + 4}{32K - 8}} \,
Q^\frac{1}{8K-2} \,
\left( q' \right)^\frac{k+1}{4K-1} \,
d(Q)^\frac{3k^2 + 11k + 8}{16K - 4} \,
(\log Q)^\frac{k + 3}{8K - 2}
$$
where $d(Q)$ is the number of divisors of $Q$, and the implicit constant is absolute.
\label{Thm5}
\end{thm}

Our proof of this is a straightforward extension of the arguments
given in \cite{G-R}. For the sake of completeness, we write out all
the necessary modifications explicitly in the appendix.\newline

Armed with Theorem \ref{Thm5}, we deduce Lemma 4 in short order.
Set
$$
T := \mathcal{P}(Q)^{3k} \, Q^\frac{1}{k} \, (q')^2 \, d(Q)^{3k} \, (\log
Q)^\frac{16 K}{k} .
$$
If $T \leq Q$, then for all $t \geq T$ Theorem \ref{Thm5} implies
$$
\left| \sum_{n \leq t} \chi(n) \right| \ll \frac{t}{\log Q}
$$
whence
$$
\left| \int_T^Q \frac{1}{t^2} \left(\sum_{n \leq t} \chi(n)\right)
\dx t \right| \ll 1.
$$
From the bound (\ref{LsChi}) we deduce that for $T \leq Q$, $|L(s_y,
\chi)| \ll \log T$. But for $T > Q$ such a bound holds trivially
(irrespective of our choice of $T$). Therefore
$$
\left| L(s_y,\chi) \right| \ll \log T \ll k \log \mathcal{P}(Q) +
\frac{1}{k}\,\log Q + \log q' + k \log d(Q) + \frac{K}{k} \, \log
\log Q .
$$
It remains to choose $k$ appropriately. Let
$$
k' := \min\left\{ \frac{1}{10} \log \log Q, \sqrt{\frac{\log Q}{\log
\mathcal{P}(Q) + \log d(Q)}} \, \right\} ,
$$
and set $k = [k'] + 1$. Writing $K' = 2^{k'}$ we have
$$
k' \log \mathcal{P}(Q) + k' \log d(Q) \ll
\sqrt{(\log Q) \big(\log \mathcal{P}(Q) + \log d(Q)\big)} \ll
\frac{1}{k'} \, \log Q
$$
and
$$
\frac{K'}{k'} \, \log \log Q \ll
(\log Q)^\frac{\log 2}{10} (\log \log Q) \ll
(\log Q)^\frac{1}{10} \ll
\frac{1}{k'} \, \log Q
$$
Finally, since $K \ll K'$ and $k \asymp k'$ (i.e. $k \ll k' \ll k$)
for all $Q$ sufficiently large, we deduce:
$$
\left| L(s_y, \chi) \right| \ll \log q' + \frac{1}{k} \, \log Q \ll
\log q' + \frac{\log Q}{\log \log Q} + \sqrt{(\log Q)\big(\log \mathcal{P}(Q) + \log
d(Q)\big)} .
$$
The proof of Lemma 4 is now complete.

\section{Proof of Lemma 5} \label{powerfullL}

Iwaniec, inspired by work of Postnikov \cite{Pos} and Gallagher
\cite{Gal}, proved the following:
\begin{thm}[See Lemma 6 of \cite{Iw}]
Given $\chi$ (mod $Q$) a primitive Dirichlet character. Then for all
$N,N'$ satisfying $(\text{rad } Q)^{100} < N < 9Q^2$ and $N < N' <
2N$,
$$
\left| \sum_{N \leq n \leq N'} \chi(n) \right| < \gamma_{_N} \, N^{1
- \epsilon_{_N}}
$$
where
$$
\gamma_x := \exp(C_1 z_x \log^2 C_2 z_x) \qquad \qquad \epsilon_x :=
\frac{1}{C_3 z_x^2 \log C_4 z_x} \qquad \qquad z_x := \frac{\log
3Q}{\log x}
$$
and the $C_i$ are effective positive constants independent of $Q$.
\end{thm}

In fact, Lemma 6 of \cite{Iw} is more general (bounding sums of 
$\chi(n) \, n^{it}$), and provides explicit choices of the 
constants $C_i$.

\begin{proof}[Proof of Lemma 5:]$\\$
Recall the bound (\ref{LsChi}):
$$
\left| L(s_y,\chi) \right| \ll \left| \int_T^Q \frac{1}{t^2}
\left(\sum_{n \leq t} \chi(n)\right) \dx t \right| + 1 + \log T .
$$
Writing
$$
\left| \sum_{n \leq t} \chi(n) \right| \leq \sqrt{t} +
\left|\sum_{\sqrt{t} < n \leq t} \chi(n)\right| ,
$$
partitioning the latter sum into dyadic intervals, and applying
Iwaniec's result to each of these, we deduce that so long as
$\sqrt{t} > (\text{rad } Q)^{100}$,
$$
\left|\sum_{n \leq t} \chi(n)\right| \ll (\log t) \, \gamma_t \,
t^{1 - \epsilon_t}
$$
with $C_1 = 400,\, C_2 = 2400,\, C_3 = 4 \cdot 1800^2,\, C_4 = 7200$
in the definitions of $\gamma_t$ and $\epsilon_t$. Choosing $T =
\exp((\log Q)^\alpha)$ for some $\alpha \in (0,1)$ to be determined
later, and assuming that $T > (\text{rad } Q)^{200}$, our bound
becomes
\begin{equation}
\left| L(s_y,\chi) \right| \ll (\log Q)^\alpha + \int_{\exp((\log
Q)^\alpha)}^Q \frac{\log t}{t^2} \, \gamma_t \, t^{1-\epsilon_t} \dx
t \label{withalphabound}
\end{equation}
Denote by $\int$ the integral in (\ref{withalphabound}), and set
$\delta_{_Q} = \frac{\log 3}{\log Q}$. Making the substitution $z =
\frac{\log 3Q}{\log t}$ and simplifying, one finds
\begin{eqnarray*}
\int & = & (\log^2 3Q) \, \int_{1 +
\delta_{_Q}}^{(1+\delta_{_Q})(\log Q)^{1-\alpha}} \frac{1}{z^3} \,
\exp\left(C_1 z \log^2 C_2 z - \frac{\log 3Q}{C_3 z^3 \log C_4 z}
\right) \dx z \\
& \ll & \exp\left( 2\log\log 3Q + C_1 (\log Q)^{1-\alpha} (\log \log
Q)^2 - \frac{(\log Q)^{3\alpha - 2}}{C_3 \log \log Q} \right)
\int_{1 + \delta_{_Q}}^{(1+\delta_{_Q})(\log Q)^{1-\alpha}}
\frac{\dx z}{z^3} \\
& \ll & 1
\end{eqnarray*}
upon choosing $\alpha = \frac{3}{4} + \epsilon$. Plugging this back
into (\ref{withalphabound}) we conclude.
\end{proof}
It is plausible that with a more refined upper bound on the integral
in (\ref{withalphabound}) one could take a smaller value of
$\alpha$, thus improving the exponents in both Lemma 5 and Theorem
2.

\section{Upper bounds on character sums} \label{upperbound}

Given \chimodq a primitive character, recall from Proposition \ref{TheBound} the 
bound
$$
\left|S_\chi(x)\right|
\ll
\sqrt{q} \, \sqrt{(\log q) \, \left|L(s_q,\chi \, \overline{\xi})\right|} +
\sqrt{q} \, (\log q)^{6/7}
$$
where $\xi$ (mod $m$) is the primitive character with $m < (\log
q)^{1/3}$ which $\chi$ is closest to, and $s_q := 1 + \frac{1}{\log q}$ . \newline

To prove Theorems 1 and 2, we would like to apply Lemmas 4 and 5 (respectively)
to derive a bound on $|L(s_q, \chi \overline{\xi})|$. An immediate difficulty is 
that both lemmas require the character to be primitive, which is not necessarily
true of $\chi \overline{\xi}$. Instead, we will apply the lemmas to the
primitive character which induces $\chi \overline{\xi}$; thus, we must understand
the size of the conductor of $\chi \overline{\xi}$. This is the goal of the
following simple lemma, which is surely well-known to the experts but which the 
author could not find in the literature. We write $[a,b]$ to denote the least 
common multiple of $a$ and $b$, and cond$(\psi)$ to denote the conductor of a 
character $\psi$.

\begin{lem}
\label{condlem}
For any non-principal Dirichlet characters $\chi_1 \, (\text{mod } q_1)$ and 
$\chi_2 \, (\text{mod } q_2)$,
$$
\text{cond}(\chi_1 \chi_2) \, \Big\vert \, \big[ \text{cond}(\chi_1) , \text{cond}(\chi_2) \big]
$$
\end{lem}

\begin{proof}
First, observe that $\chi_1 \chi_2$ is a character modulo $[q_1, q_2]$: one needs only check
that it is completely multiplicative, periodic with period $[q_1, q_2]$, and that
$\chi_1 \chi_2 (n) = 0$ if and only if $(n,[q_1, q_2]) > 1$. Since the conductor of a 
character divides its modulus, the lemma is proved in the case that both $\chi_1$ and
$\chi_2$ are primitive. \newline

Now suppose that $\chi_1$ and $\chi_2$ are not necessarily primitive; denote by
$\tilde{\chi}_i \, (\text{mod } \tilde{q}_i)$ the primitive character which induces $\chi_i$.
By the argument above, we know that
\begin{equation}
\text{cond}(\tilde{\chi}_1 \tilde{\chi}_2) \, \Big\vert \, \big[\tilde{q}_1 , \tilde{q}_2 \big] .
\label{tilde}
\end{equation}
Next we note that the character $\tilde{\chi}_1 \tilde{\chi}_2$, while not necessarily primitive,
does induce $\chi_1 \chi_2$ (i.e. $\chi_1 \chi_2 = \tilde{\chi}_1 \tilde{\chi}_2 \chi_0$ for
$\chi_0$ the trivial character modulo $[q_1, q_2]$), whence 
$\text{cond}(\tilde{\chi}_1 \tilde{\chi}_2) = \text{cond}(\chi_1 \chi_2)$.
Plugging this into (\ref{tilde}) we immediately deduce the lemma.
\end{proof}

Given \chimodq and $\xi \, (\text{mod } m)$ as at the start of the section, denote by 
$\psi \, (\text{mod } Q)$ the primitive character inducing $\chi \overline{\xi}$. Taking 
$\chi_1 = \chi$ and $\chi_2 = \overline{\xi}$ in Lemma \ref{condlem}, we see that
$Q \mid [q,m]$; in particular, $Q \leq qm$. On the other hand, making the choice 
$\chi_1 = \chi \overline{\xi}$ and $\chi_2 = \xi$ yields $q \mid [Q,m]$, so $q \leq Qm$.
Combining these two estimates, we conclude that
\begin{equation}
\label{condbound}
\frac{q}{m} \leq Q \leq qm
\end{equation}

Since we will be working with both $L(s, \chi \overline{\xi})$ and $L(s, \psi)$, the
following estimate will be useful:
\begin{lem}
\label{Lcompare}
Given \chimodq and $\xi \, (\text{mod } m)$ primitive characters, let 
$\psi \, (\text{mod } Q)$ be the primitive character which induces $\chi \overline{\xi}$.
Then for all $s$ with Re$(s) > 1$,
$$
\left| \frac{L(s, \chi \overline{\xi})}{L(s, \psi)} \right| \ll
1 + \log \log m .
$$
\end{lem}

\begin{proof}
For Re$(s) > 1$ we have
$$
\frac{L(s, \chi \overline{\xi})}{L(s, \psi)} = 
\mathop{\prod_{p \mid [q,m]}}_{p \nmid Q} \left( 1 - \frac{\psi(p)}{p^s} \right)
$$
whence
$$
\left| \frac{L(s, \chi \overline{\xi})}{L(s, \psi)} \right| \leq
\mathop{\prod_{p \mid [q,m]}}_{p \nmid Q} \left( 1 + \frac{1}{p} \right) .
$$
From Lemma \ref{condlem} we know $q \mid [Q,m]$; it follows that if 
$p \mid [q,m]$ and $p \nmid Q$ then $p$ must divide $m$. Thus,
$$
\mathop{\prod_{p \mid [q,m]}}_{p \nmid Q} \left( 1 + \frac{1}{p} \right)
\leq \prod_{p \mid m} \left( 1 + \frac{1}{p} \right) .
$$
Since
$$
\log \, \prod_{p \mid m} \left( 1 + \frac{1}{p} \right) =
\sum_{p \mid m} \log \left( 1 + \frac{1}{p} \right) \leq
\sum_{p \mid m} \frac{1}{p} ,
$$
to prove the lemma it suffices to show that for all $m$ sufficiently large,
\begin{equation}
\sum_{p \mid m} \frac{1}{p} \leq \log \log \log m + O(1) .
\label{sumprimediv}
\end{equation}
Let $P = P(m)$ denote the largest prime such that $\prod_{p \leq P} p \leq m$.
Then $\omega(m) \leq \pi(P)$ (otherwise we would have 
$m \geq \text{rad}(m) > \prod_{p \leq P} p$, contradicting the maximality of $P$); 
therefore,
$$
\sum_{p \mid m} \frac{1}{p} \leq \sum_{p \leq P} \frac{1}{p} = \log \log P + O(1) .
$$
Finally from the prime number theorem we know that for all $m$ sufficiently large,
$\theta(P) \geq \frac{1}{2} P$, whence $P \leq 2 \log m$ and the bound
(\ref{sumprimediv}) follows.
\end{proof}

With these lemmas in hand we can now prove Theorems 1 and 2 without too much difficulty.

\begin{proof}[Proof of Theorem 1:]
Given \chimodq primitive with $q$ squarefree, define the character $\xi \, (\text{mod } m)$
as in Theorem A, and let $\psi \, (\text{mod } Q)$ be the primitive character inducing
$\chi \overline{\xi}$. Recall that we denote the largest prime factor of $n$ by
$\mathcal{P}(n)$.\newline

From Proposition \ref{TheBound} we have 
\begin{equation}
\label{thm1stp1}
\left|S_\chi(x)\right|
\ll
\sqrt{q} \, \sqrt{(\log q) \, \left|L(s_q,\chi \, \overline{\xi})\right|} +
\sqrt{q} \, (\log q)^{6/7} ,
\end{equation}
and Lemma \ref{Lcompare} yields the bound
\begin{equation}
\label{thm1stp2}
\left|L(s_q,\chi \, \overline{\xi})\right| \ll |L(s_q,\psi)| \, \log \log \log q .
\end{equation}
Lemma \ref{condlem} tells us that $Q \mid [q,m]$, whence for all primes $p > m$ we have
$$
\text{ord}_p \, Q \leq \max (\text{ord}_p \, q , \text{ord}_p \, m) = \text{ord}_p \, q \leq 1
$$
since $q$ is squarefree. Therefore we may apply Lemma 4 to the character $\psi$, taking $y = q$ and
$$
q' = \prod_{p \leq m} p^{\text{ord}_p Q} \, ;
$$
this gives the bound
$$
|L(s_q,\psi)| \ll \log q' + \frac{\log Q}{\log \log Q} + 
\sqrt{(\log Q) \log \big(\mathcal{P}(Q) d(Q) \big)} .
$$
It remains only to bound the right hand side in terms of $q$, which we do term by term.
The first term is small:
\begin{eqnarray*}
\log q' 
& = & 
\sum_{p \leq m} (\text{ord}_p \, Q) \log p \\
& \leq &
\sum_{p \leq m} (\text{ord}_p \, q) \log p + \sum_{p \leq m} (\text{ord}_p \, m) \log p \\
& \leq &
\theta(m) + \log m \\
& \ll & 
(\log q)^\frac{1}{3} .
\end{eqnarray*}
From (\ref{condbound}) we deduce
$$
\frac{\log Q}{\log \log Q} \ll \frac{\log q}{\log \log q} \, .
$$
For the last term, Lemma \ref{condlem} yields
$$
d(Q) \leq d(qm) \leq d(q) d(m) \leq d(q) (\log q)^\frac{1}{3}
$$
and
$$
\mathcal{P}(Q) \leq \max \big( \mathcal{P}(q), \mathcal{P}(m) \big)
\leq \mathcal{P}(q) \mathcal{P}(m) \leq \mathcal{P}(q) (\log q)^\frac{1}{3} \, ,
$$
while (\ref{condbound}) gives $\log Q \ll \log q$.
Putting this all together, we find
$$
|L(s_q,\psi)| \ll \frac{\log q}{\log \log q} + 
\sqrt{(\log q) \log \big(\mathcal{P}(q) d(q) \big)} \, ;
$$
plugging this into (\ref{thm1stp2}) and (\ref{thm1stp1}) we deduce the theorem.
\end{proof}

\begin{proof}[Proof of Theorem 2:]
Given \chimodq with $q$ large and rad$(q) \leq \exp\left((\log q)^\frac{3}{4}\right)$, 
let $\xi \, (\text{mod } m)$ be defined
as in Theorem A, and let $\psi \, (\text{mod } Q)$ denote the primitive 
character which induces $\chi \overline{\xi}$. We have rad$(m) \leq
\exp\big(\theta(m)\big)$, whence by the prime number
theorem $\exists \, C > 0$ with
\begin{eqnarray*}
\text{rad}(Q) & \leq & \text{rad}(q) \, \text{rad}(m) \\
& \leq & \exp\left((\log q)^{3/4} + C \, (\log q)^{1/3}\right) \\
& \leq & \exp\left(\frac{4}{3} \, (\log q)^{3/4}\right) 
\end{eqnarray*}
for all $q$ sufficiently large. From (\ref{condbound}) we deduce
$$
\left(\frac{\log Q}{\log q}\right)^\frac{3}{4}
\geq 
\left(\frac{\log \frac{q}{m}}{\log q}\right)^\frac{3}{4}
\geq
\left(1 - \frac{\log \log q}{\log q}\right)
\geq
\frac{2}{3}
$$
for $q$ sufficiently large, whence
$$
\text{rad}(Q) \leq \exp\left(2 \, (\log Q)^{3/4}\right) \, .
$$
Combining Lemma 5 with (\ref{thm1stp2}) and (\ref{condbound}) we obtain
$$
\left| L(s_q,\, \chi \, \overline{\xi}) \right| 
\ll_\epsilon (\log \log \log q) (\log Q)^{3/4 + \epsilon} 
\leq (\log \log \log q) (\log qm)^{3/4 + \epsilon} 
\ll_\epsilon (\log q)^{3/4 + \epsilon} .
$$
Plugging this into Proposition \ref{TheBound} yields Theorem 2.
\end{proof}

\appendix

\section{Appendix: Proof of Theorem \ref{Thm5}}

We follow the original proof of Graham and Ringrose very closely;
indeed, we will only explicitly write down those parts of their
arguments which must be modified to obtain our version of the
result. We refer the reader to sections 3 - 5 of \cite{G-R}. Set $S
:= \sum_{M < n \leq M+N} \chi(n)$. \newline

We begin by restating Lemma 3.1 of \cite{G-R}, but skimming off some
of the unnecessary hypotheses given there:
\begin{lem}[Compare to Lemma 3.1 of \cite{G-R}]
Let $k \geq 0$ be an integer, and set $K := 2^k$. Let
$q_0,\ldots,\,q_k$ be arbitrary positive integers, and let $H_i :=
N/q_i$ for all $i$. Then
\begin{equation}
\left|S\right|^{2K} \leq 8^{2K-1} \left( \max_{0 \leq j \leq k}
\left(N^{2K - K/J} \, q_j^{K/J}\right) + \frac{N^{2K-1}}{H_0 \cdots
H_k} \sum_{h_0 \leq H_0} \cdots \sum_{h_k \leq H_k} \left| S_k({\bf
h}) \right| \right) \label{lemma3.1}
\end{equation}
where $J = 2^j$ and $S_k({\bf h})$ satisfies the bound given below.
\label{3.1'}
\end{lem}
A bound on $S_k({\bf h})$ is given by (3.4) of \cite{G-R}:
\begin{equation}
|S_k({\bf h})| \ll N \, Q^{-1} \, |S(Q;\, \chi,\, f_k,\, g_k,\, 0)|
+ \sum_{0 < |s| \leq Q/2} \frac{1}{|s|} \, \left| S(Q; \, \chi,\,
f_k,\, g_k,\, s) \right| . \label{3.4'}
\end{equation}
See pages 279-280 of \cite{G-R} for the definitions of $f_k,g_k$,
and $S(Q; \chi, f_k, g_k, s)$.
\newline

Let $q := Q/q'$. We have $(q,q') = 1$, whence from Lemma 4.1 of
\cite{G-R} we deduce
$$
S(Q; \chi, f_k, g_k, s) = S(q'; \chi', f_k, g_k, s \overline{q})
S(q; \eta, f_k, g_k, s \overline{q'})
$$
for some primitive characters $\chi'$ (mod $q'$) and $\eta$ (mod
$q$), where $q \overline{q} \equiv 1$ (mod $q'$) and $q'
\overline{q'} \equiv 1$ (mod $q$). By construction, $q$ is
squarefree, so Lemmas 4.1-4.3 of \cite{G-R} apply to give
$$
\left|S(q;\,\eta,\,f_k,\,g_k,\,s \overline{q'})\right| \leq
d(q)^{k+1} \, \left(\frac{q}{(q,Q_k)}\right)^{1/2} \, (q,Q_k,|s
\overline{q'}|)
$$
where $Q_k := \prod_{i \leq k} h_i q_i$. Combining this with the
trivial estimate $\left| S(q'; \chi', f_k, g_k, s \overline{q})
\right| \leq q'$ yields:
\begin{lem}[Compare with Lemma 4.4 of \cite{G-R}]
Keep the notation as above. Then for any positive integers
$q_1,\ldots,\,q_k$,
$$
\left|S(Q;\,\chi,\,f_k,\,g_k,\,s)\right| \leq q' \, d(q)^{k+1} \,
\left(\frac{q}{(q,Q_k)}\right)^{1/2} \, (q,Q_k,|s \overline{q'}|) .
$$
\end{lem}

We shall need the following simple lemma (versions of which appear
implicitly in \cite{G-R}):
\begin{lem}
\label{4.4}
Given $q$, $\overline{q'}$ be as above; let $x$ and $H$ be
arbitrary. Then
\begin{eqnarray}
&1.& \label{lempart1} 
\sum_{0 < |s| \leq x} \frac{(q,\, |s \overline{q'}|)}{|s|} \ll 
d(q) \, \log x  \\
&2.& \label{lempart2}
\sum_{h \leq H} (q,\, h)^\frac{1}{2} \leq d(q) \, H
\end{eqnarray}
\end{lem}
\begin{proof} $\\$
\begin{enumerate}
\item Since $(q,q')=1$, we have $(q,\overline{q'})=1$, whence
$$
\sum_{0 < |s| \leq x} \frac{(q,\, |s \overline{q'}|)}{|s|} = 2 \,
\sum_{1 \leq s \leq x} \frac{(q,\, s \overline{q'})}{s} = 2 \,
\sum_{1 \leq s \leq x} \frac{(q,\, s)}{s} = 2 \, \sum_{n \geq 1}
\frac{a_n}{n}
$$
where $a_n := \#\{s \leq x : n = \frac{s}{(q,\, s)}\}$. Note that
$a_n = 0$ for all $n > x$, and that
$$
a_n = \#\{s \leq x : s = (q,\, s) \, n\} \leq \#\{s \leq x : s = d
n,\, d \mid q\} \leq d(q) .
$$
Therefore
$$
\sum_{0 < |s| \leq x} \frac{(q,\, |s \overline{q'}|)}{|s|} \ll
\sum_{n \geq 1} \frac{a_n}{n} \leq d(q) \, \sum_{n \leq x}
\frac{1}{n} \ll d(q) \log x .
$$
\item Write
$$
\sum_{h \leq H} (q,\, h)^\frac{1}{2} = \sum_{n \geq 1} a_n \,
\sqrt{n}
$$
where $a_n := \#\{h \leq H : n = (q,\, h)\}$. It is clear that $a_n
= 0$ whenever $n \nmid q$. Also, if $(q,\, h) = n$ then $n \mid h$,
whence
$$
a_n \leq \#\{h \leq H : n \mid h\} \leq \frac{H}{n} .
$$
Therefore
$$
\sum_{h \leq H} (q,\, h)^\frac{1}{2} = \sum_{n \geq 1} a_n \,
\sqrt{n} \leq \sum_{n \mid q} \frac{H}{\sqrt{n}} \leq d(q) \, H .
$$
\end{enumerate}
\end{proof}

\begin{lem}[Compare to Lemma 4.5 of \cite{G-R}]
Keep the notation from above. For any real number $A_0 \geq 1$,
$$
|S|^{4K} \ll 8^{4K-2} \left( A \, A_0^{2K} + B \, A_0^{-2K + 1} \,
(q')^2 + C \, A_0^{2K-1} \, (q')^2 \right)
$$
where
\begin{eqnarray*}
A &=& N^{2K} \\
B &=& N^{6K - k - 4} \, P^{k+1} \, Q \, d(Q)^{2k+4} \, \log^2 Q \\
C &=& N^{2K + k + 2} \, Q^{-1} \, d(Q)^{4k+4}
\end{eqnarray*}
and the implied constant is independent of $k$.
\end{lem}
(Note that in the original paper, there is a persistent typo of
writing $M$ rather than $N$.)
\begin{proof}
Following the proof of Lemma 4.5 in \cite{G-R} and applying
(\ref{lempart1}) with $x = Q/2$ yields the following analogue of
equation (4.5) from that paper:
\begin{equation}
\sum_{h_k \leq H_k} \sum_{0 < |s| \leq Q/2} \frac{1}{|s|} \, \left|
S(Q; \, \chi,\, f_k,\, g_k,\, s) \right| \ll q' \, \sqrt{q} \,
d(q)^{k+2} \, H_k \, R_k^{-\frac{1}{2}} \, \log Q . \label{4.5'}
\end{equation}
Setting $S_j := h_0 \cdots h_j$, one deduces the following analogue
of equation (4.6) of \cite{G-R}:
$$
N \, Q^{-1} \, \sum_{h_k \leq H_k} |S(Q;\, \chi,\, f_k,\, g_k,\, 0)|
\leq N \, q' \, \frac{\sqrt{q \, R_k}}{Q} \, d(Q)^{k+2} \, H_k \,
\sqrt{(q,\, S_{k-1})} .
$$
From (\ref{lempart2}) and the bound $(q,\, S_j) \leq (q,\, S_{j-1})
\, (q,\, h_j)$ one sees that
\begin{equation}
\sum_{h_0 \leq H_0} \cdots \sum_{h_{k-1} \leq H_{k-1}} \sqrt{(q,\,
S_{k-1})} \leq d(q)^k \, H_0 \cdots H_{k-1} . \label{4.6'}
\end{equation}
Plugging (\ref{3.4'}) into (\ref{lemma3.1}) and applying
(\ref{4.5'}) and (\ref{4.6'}), one obtains:
\begin{eqnarray*}
|S|^{2K} & \ll &
    8^{2K-1} \, \max_{0 \leq j \leq k} \left(N^{2K - K/J} \, q_j^{K/J}\right) + \\
& & 8^{2K-1} \, q' \, N^{2K-1} \, d(q)^{k+2} \, (\log Q) \, \sqrt{\frac{q}{R_k}} + \\
& & 8^{2K-1} \, q' \, N^{2K} \, d(q)^{2k+2} \, \frac{\sqrt{q}}{Q} \,
\sqrt{R_k} .
\end{eqnarray*}
Since $q \mid Q$, we have that $q \leq Q$ and $d(q) \leq d(Q)$.
Therefore from the above we deduce the following analogue of (4.7)
in \cite{G-R}:
\begin{eqnarray*}
|S|^{2K} & \ll &
    8^{2K-1} \, \max_{0 \leq j \leq k} \left(N^{2K - K/J} \, q_j^{K/J}\right) + \\
& & 8^{2K-1} \, q' \, N^{2K-1} \, d(Q)^{k+2} \, (\log Q) \, \sqrt{\frac{Q}{R_k}} + \\
& & 8^{2K-1} \, q' \, N^{2K} \, d(Q)^{2k+2} \, \sqrt{\frac{R_k}{Q}}
.
\end{eqnarray*}
The rest of the proof given in \cite{G-R} can now be copied exactly
to yield our claim.
\end{proof}
Chasing through the arguments in \cite{G-R} gives this analogue of
Lemma 5.3, which we record for reference:
\begin{lem}[Compare to Lemma 5.3 of \cite{G-R}]
\begin{eqnarray*}
|S| & \ll & N^{1 - \frac{k+3}{8K-2}} \, P^\frac{k+1}{8K-2} \,
Q^\frac{1}{8K-2} \, d(Q)^\frac{k+2}{4K-1} \, (\log Q)^\frac{1}{4K-1}
\, (q')^\frac{1}{4K-1}
+ \\
& & N^{1 - \frac{1}{4K}} \, P^\frac{k+1}{8K} \, d(Q)^\frac{3k+4}{4K}
\, (\log Q)^\frac{1}{4K} \, (q')^\frac{1}{2K} .
\end{eqnarray*}
\end{lem}
Finally, we arrive at:
\begin{proof}[Proof of Theorem \ref{Thm5}]
Let $E_k$ be the right hand side of the bound claimed in the
statement of the theorem. The rest of the proof given in \cite{G-R}
now goes through almost verbatim.
\end{proof}

This concludes the proof of Theorem \ref{Thm5}. Note that one can
extend this to a bound on all non-principal characters by following
the argument given directly after Lemma 5.4 in \cite{G-R}; however,
for our applications the narrower result suffices.

\textsc{Department of Mathematics, University of Michigan, Ann Arbor, MI 48109}
\newline
\textit{Email:} \texttt{lgoldmak@umich.edu}
\newline

\textit{Current address:} Department of Mathematics, Stanford
University, Bldg 380, 450 Serra Mall, Stanford, CA  94305

\end{document}